\documentclass[11pt, reqno]{article}
\usepackage{amsfonts}
\usepackage{amsmath}
\usepackage{amsthm}
\usepackage{amssymb}
\usepackage{oldgerm}
\usepackage{lastpage}
\usepackage{fancyhdr}
\usepackage{fancybox}
\usepackage{mathrsfs}
\usepackage{epsfig}
\def\f{\mathfrak }
\def\b{\mathbb }
\def\phi{\varphi }

\swapnumbers
\theoremstyle{plain}
\newtheorem{theorem}{Theorem}[section]
\newtheorem{corollary}[theorem]{Corollary}
\newtheorem{lemma}[theorem]{Lemma}
\newtheorem{proposition}[theorem]{Proposition}
\theoremstyle{definition}

\theoremstyle{remark}
\newtheorem*{remark}{Remark}

\numberwithin{equation}{section}

\begin{document}
\title{Positivity of Dunkl's intertwining operator via the 
trigonometric setting}
\author{Margit R\"osler  and Michael Voit}
\date{}

\maketitle

\begin{abstract}
In this note, a new proof for the positivity of Dunkl's intertwining operator in the crystallographic case is given. It is based on  an asymptotic relationship between the Opdam-Cherednik kernel and the Dunkl kernel as recently observed by M. de Jeu, and  on positivity results of S. Sahi for the Heckman-Opdam polynomials and their non-symmetric counterparts. 
\end{abstract}


\section{Introduction}

In \cite{R1}, it was proven that Dunkl's intertwining operator between the rational Dunkl operators for a fixed finite reflection group and nonnegative multiplicity function is positive. As a consequence we obtained an abstract Harish-Chandra type integral representation for the Dunkl kernel, the image of the usual exponential kernel under the intertwiner. 
The proof was based on methods from the theory of operator semigroups and 
 a rank-one reduction.

In the present note, we give a new, completely different proof of these
results under the only additional assumption that the underlying reflection group has to be crystallographic. In contrast to the proof of \cite{R1}, where 
precise information on the supports of the representing measures  could only be obtained by going back to estimates of the kernel from \cite{dJ1}, this
information is now directly obtained. Our new approach relies first
on an asymptotic relationship between the Opdam-Cherednik kernel and the Dunkl kernel as recently observed by de Jeu \cite{dJ2}, and second on positivity results of Sahi \cite{S} for the Heckman-Opdam polynomials and their non-symmetric counterparts. 

\section{Preliminaries}
\subsection{Basic notation}
Let $\f a $ be a finite-dimensional Euclidean 
vector space with inner product 
$\langle\,.\,,\,.\,\rangle$. We use the same notation for the bilinear
extension of $\langle\,.\,,\,.\,\rangle$ to the complexification 
$\f a_{\b C}$ of $\f a$,  and we  identify $\f a$   with its dual 
$\f a^* = \text{Hom} (\f a, \b R)$ via the given inner product. For $\alpha\in \f a\setminus\{0\}$ we write
$ \alpha^\vee = \frac{2\alpha}{\langle\alpha,\alpha\rangle}\,$
 and 
$\, \sigma_\alpha(x) = x- \langle x,\alpha^\vee\rangle
\alpha\,$
for the orthogonal reflection in the hyperplane perpendicular to $\alpha$. 
We  consider a crystallographic root system $R$ in $\f a$, i.e. $R$ is a finite subset of $\f a \setminus\{0\}$ which spans $\f a$ and satisfies
$\sigma_\alpha(R) = R$ and $\langle\alpha,\beta\rangle\in \b Z$ for all $\alpha,\beta\in R$. We
also assume that $R$ is indecomposable and reduced, i.e.  $ R\cap \b R\alpha = \{\pm\alpha\}$ for all $\alpha\in R$. Let $W$ be the finite reflection group generated by the $\sigma_\alpha,\, \alpha\in R$.
We shall fix 
a positive subsystem $R_+$ of $R$ as well as a 
\emph{nonnegative} multiplicity function $k = (k_\alpha)_{\alpha\in R}$, 
satisfying $k_\alpha =k_\beta$ if $\alpha$ and $\beta$ are in the same 
$W$-orbit.

\subsection{Rational Dunkl operators and Dunkl's intertwiner}
References for this section are \cite{D1}, \cite{D2}, \cite{O3}, \cite{R1}.
Let $\mathcal P = \b C[\f a]$ denote the vector space of complex 
polynomial functions on $\f a$, and $\mathcal P_n\subset P$ the subspace 
of polynomials which are homogeneous of degree $n\in \b Z_+$.
The rational Dunkl operators on $\f a$ associated with  $R$ and 
fixed multiplicity  $k\geq 0$ are given by
\[ T_\xi = T_\xi(k) = \partial_\xi + \sum_{\alpha\in R_+} k_\alpha 
\langle\alpha,\xi\rangle \frac{1}{\langle\alpha, \,.\,\rangle}(1-\sigma_\alpha), \quad\xi\in \f a.\] 
These operators commute and map $\mathcal P$ onto itself. Moreover, there 
exists a unique linear isomorphism $V= V_k$ on $\mathcal P$  with $V(1) = 1,\,
V(\mathcal P_n) = P_n$  and  $T_\xi V = V \partial _\xi$ for all $\xi\in \f a$.  According to \cite{D2}, the intertwining operator $V$ can be extended to larger classes of functions, as follows: For $r>0$, let 
$K_r:= \{x\in \f a: |x|\leq r \}$ denote the ball of radius $r$, and define 
\[ A_r:= \{ f: K_r\to \b C, \,f=\sum_{n=0}^\infty f_n \,\,\text{ with }\, f_n \in \mathcal P_n \text{ and }\,\, \|f\|_{A_r}:= \sum_{n=0}^\infty \|f_n\|_{\infty, K_r} < \infty\},\]
where $\|f_n\|_{\infty, K_r}:= \sup_{x\in K_r}|f_n(x)|.$
The space $A_r$ is a Banach space with norm $\|\,.\,\|_{A_r}$
(in fact, a commutative Banach algebra).
$V$ extends to a continuous linear operator on $A_r$ by $\, V\bigl(\sum_{n=0}^\infty f_n\bigr):= \sum_{n=0}^\infty Vf_n\,.$ The Dunkl kernel
 $\text{Exp}_W$ is defined by
\[\text{Exp}_W(\,.\,,z) := V(e^{\langle\,.\,,z\rangle}), \quad z\in \f a_{\b C}.\]
It extends to a holomorphic function on $\f a_{\b C}\times \f a_{\b C}$ 
which is symmetric in its arguments. For  $\lambda\in \f a_{\b C}, \,\text{Exp}_W(\lambda,\,.\,)$ is the unique holomorphic solution of the joint eigenvalue problem 
\[ T_\xi f\,=\,\langle\lambda,\xi\rangle f \quad\forall \xi\in \f a, \quad f(0)=1. \]
For $x\in\f a$, we denote by $C(x)$ the closure of the 
convex hull of the $W$-orbit $Wx$ of $x$ in $\f a$. Moreover, for a locally compact Hausdorff space $X$, we write $M^1(X)$ for the set of probability measures 
on the Borel-$\sigma$-algebra of $X$. 
In \cite{R1}, the following is proven:

\begin{theorem}\label{postheorem} For each $x\in \f a$ there exists a unique  probability
measure $\mu_{x}\in M^1(\f a)$ such that
\[ Vf(x) = \int_{\f a} f(\xi) d\mu_x(\xi) \quad\forall\, f\in A_{|x|}.\]
The support of $\mu_x$ is contained in $C(x)$.
\end{theorem}

As a consequence, 
\begin{equation}\label{kernelint}
\text{Exp}_W(x,z) = \int_{\f a} e^{\langle\xi,z\rangle} 
 d\mu_x(\xi) \quad \forall z\in \f a_{\b C}.
\end{equation}
In \cite{R1}, the proof of the inclusion  $\text{supp}\,\mu_x \subseteq C(x)$
 requires the exponential  bounds on 
$\text{Exp}_W$ from \cite{dJ1}, which are by far not straightforward. As is well-known, the $W$-invariant parts of the rational and trigonometric Dunkl theories are, for certain discrete sets of multiplicities, realized within  the classical Harish-Chandra theory
for semisimple symmetric spaces. In particular, for such $k$ the generalized 
Bessel functions 
\[ J_W(\,.\,,z) = \frac{1}{|W|} \sum_{w\in W} \text{Exp}_W(\,.\,, w^{-1}z), \quad z\in \f a_{\b C}\]
 can be identified with the spherical functions  of 
an underlying Cartan motion group, for details see e.g. \cite{dJ2}. 
In this case, their integral representation according to  \eqref{kernelint} 
is a special case of a (Euclidean type) Harish-Chandra integral,  
 and the inclusion $\text{supp}\,\mu_x \subseteq C(x)$ 
follows from Kostant's convexity theorem (\cite{He},  Propos. 
IV.4.8 and Theorem IV.10.2).

\subsection{Cherednik operators and the Opdam-Cherednik  kernel} 

The basic concepts of this as well as the following section are due to Opdam \cite{O1} (see also Part I of \cite{O2}), Heckman and Cherednik \cite{C}. 
Let \[ P :=\{\lambda\in \f a: \langle \lambda, 
   \alpha^\vee \rangle\,\in \b Z \quad\text{for all }\,  \alpha\in R\}\]
denote the weight lattice associated with the root system $R$. 
For $\lambda \in \f a_{\b C}$ we define the exponential $e^\lambda$ on 
$\f a_{\b C}$ by  
$\, e^\lambda(z):= e^{\langle\lambda,z\rangle}$, and denote by
$\mathcal T$ the $\b C$-span of   $\{e^\lambda, \, \lambda\in P\}$. This is 
 the algebra
of trigonometric polynomials on $\f a_\b C$ with respect to $R$. 
The Cherednik operator in direction $\xi\in \f a$ is
defined by
\[D_\xi = D_{\xi}(k) = \partial_{\xi} + \sum_{\alpha\in R_+} k_\alpha \langle \alpha,\xi\rangle \frac{1}{1-e^{-\alpha}}(1-\sigma_\alpha) - 
\langle \rho(k),\xi\rangle,\]
where $\rho(k) = \frac{1}{2} \sum_{\alpha\in R_+} k_\alpha \alpha$.
Each $D_\xi$ maps $\mathcal T$ onto itself and (for
 fixed $k$) the operators $D_\xi$ commute. Notice that in contrast to the rational $T_\xi$, they depend on 
the particular choice of $R_+$. For each $\lambda\in \f a_{\b C}$, there 
exists a unique holomorphic function $G(\lambda,\,.\,)$ in a 
tubular neighborhood of $\f a$ which satisfies
\begin{equation}\label{G_char}
 D_\xi G(\lambda,\,.\,) = \langle\lambda,\xi\rangle\, G(\lambda,\,.\,)\quad 
\forall\, \xi\in \f a,  \quad G(\lambda, 0) =1,
\end{equation}
see Corollary I.7.6 of \cite{O2}. $G$ is called the Opdam-Cherednik kernel. It is in fact (as a function of both arguments) holomorphic in a suitable tubular neighborhood of 
$\f a_{\b C} \times \f a$ (\cite{O1}, Theorem 3.15).
 The rational Dunkl operators can be considered  a scaling limit 
of the Cherednik operators, and this implies limit relations between 
the kernels $\text{Exp}_W$ and $G$. We shall need the following variant of 
Theorem 4.12 in \cite{dJ2}:

\begin{proposition}\label{scalinglimit} Let $\delta >0$ be a constant, $K, L \subset \f a_{\b C}$  compact sets
and $h: (0,\delta)\times L \to \f a_{\b C}$ a continuous mapping such that 
$\,\lim_{\epsilon\to 0} \epsilon h(\epsilon,\lambda) = \lambda $ 
uniformly on $L$.
Then
\[ \lim_{\epsilon\to 0} G(h(\epsilon,\lambda), \epsilon z) \,=\, 
{\rm Exp}_W(\lambda,z)\]
uniformly for $(\lambda,z)\in L\times K.$ 
\end{proposition}

\noindent
The proof  is  the same as for Theorem 4.12 in \cite{dJ2},   with $\lambda/\epsilon$ replaced by $h(\epsilon,\lambda)$. We mention that for integral 
$k$, such a limit transition has first been carried out in \cite{BO},
by use of shift operator methods.

\subsection{A scaling limit for non-symmetric Heckman-Opdam polynomials}
The definition of these  polynomials involves a suitable partial order on $P$;
we refer to the one used in  \cite{O2}.
Let
\[ P_+ :=\{ \lambda\in P: \langle \lambda, \alpha^\vee\rangle \,\geq 0 \quad\text{for all }\,\alpha\in R_+\}\]
denote the set of dominant weights associated with $R_+$, and 
 $\lambda_+$  the unique 
 dominant weight in the orbit $W\lambda$. One defines
$\lambda  \vartriangleleft  \nu$ if either $\lambda_+ <\nu_+$ in dominance ordering 
(i.e. $\nu_+-\lambda_+ \in Q_+,$ the $\b Z_+$-span of $R_+$) or if $\lambda_+ =
\nu_+$ and $\nu <\lambda$ (in dominance ordering). Further, $\lambda \trianglelefteq \nu$ means $\lambda = \nu$
or $\lambda\vartriangleleft\nu$.
The non-symmetric Heckman-Opdam polynomials 
 $\{ E_\lambda\,: \,\lambda\in P\} \subset \mathcal T$ associated with $R_+$ and $k$ are 
uniquely  characterized by the conditions
\begin{align}
1.&\quad  E_\lambda = \sum_{\nu\trianglelefteq \lambda} a_{\lambda,\nu}\, e^\nu \,\text{ with }\,\,a_{\lambda,\lambda} = 1 \label{coeff1}\\
2.&\quad  D_\xi E_\lambda = \langle\,\widetilde\lambda,\xi\rangle\, 
E_\lambda \quad \forall \, \xi\in \f a,\label{coeff2}
\end{align}
with the shifted spectral variable 
$\, \widetilde\lambda =  \lambda + 
\frac{1}{2} \sum_{\alpha\in R_+} k_\alpha \epsilon(\langle\lambda,\alpha^\vee\rangle) \alpha.\,$
 Here $\epsilon: \b R\to \{\pm 1\}$ is defined by
$\epsilon(x) = 1$ for $x>0$, $\epsilon(x) = -1$ for $x\leq 0$. For details
see \cite{O1} and \cite{O2}, Section I.2.3.

 On the other hand, we 
know (c.f. \eqref{G_char}) that $G(\widetilde \lambda,\,.\,)$ is the up to a constant factor unique holomorphic solution of \eqref{coeff2}. Hence
\begin{equation}\label{identif}
 E_\lambda = \, c_\lambda\cdot G(\widetilde\lambda,\,.\,) 
\end{equation}
with a constant $c_\lambda = E_\lambda(0) >0$. The precise value of $c_\lambda$ is given in Theorem 4.7 of \cite{O2}.

\begin{corollary} \label{limit}
For $\lambda\in P$ and $z\in \f a_{\b C}$, 
\[{\rm Exp}_W(\lambda,z) = \lim_{n\to\infty} \frac{1}{c_{n\lambda}} E_{n\lambda}\bigl(\frac{z}{n}\bigr).\]
The convergence is locally uniform with respect to $z$.
\end{corollary}

\begin{proof}
Fix $\lambda\in P$ and observe that  
$\,\widetilde{n\lambda} = n\lambda + \frac{1}{2} \sum_{\alpha\in R_+} k_\alpha \epsilon(\langle\lambda,\alpha^\vee\rangle) \alpha\,$
 for all $n\in \b N$. Thus
by Proposition \ref{scalinglimit} and identity \eqref{identif} we have, locally uniformly for $z\in \f a_{\b C}$, 
\[ {\rm Exp}_W(\lambda,z) =\, \lim_{n\to\infty} G\bigl(\widetilde{n\lambda},\frac{z}{n}\bigr) = \, \lim_{n\to\infty} \frac{1}{c_{n\lambda}} E_{n\lambda}\bigl(\frac{z}{n}\bigr).\]
\end{proof}

\begin{remark}
Similar results hold for
 the symmetric Heckman-Opdam polynomials 
\[ P_\lambda(z) = \frac{|W\!\lambda|}{|W|} \sum_{w\in W} E_\lambda (w^{-1}z),
\quad \lambda\in P_+\,.\]
They are $W$-invariant and related with the multivariable hypergeometric function  
\[ F(\lambda,z) = \frac{1}{|W|} \sum_{w\in W} G(\lambda, w^{-1}z)\]
via 
\[ P_\lambda  = c_\lambda^*\cdot F(\lambda + \rho,\,.\,)\quad \forall
\lambda\in P_+\]
with $\,c_\lambda^* = |W\!\lambda|\cdot c_\lambda\,$ and 
 $\,\rho = \rho(k) = \frac{1}{2}\sum_{\alpha\in R_+} k_\alpha \alpha\,$, c.f. Eq.~(4.4.10) in \cite{HS}. This also follows from \eqref{identif}, because
$F$ is in fact $W$-invariant in both arguments and for $\lambda\in P_+$ the shifted weight $\widetilde\lambda$ is contained in the $W$-orbit of $\lambda +\rho$ (Proposition 2.10 of \cite{O1}).
Further, Corollary \ref{limit} implies that for $\lambda\in P_+$ and $z\in \f a_{\b C}$, 
\begin{equation}\label{limit2} J_W(\lambda,z) = \,\lim_{n\to\infty} F\bigl(n\lambda + \rho,\frac{z}{n}\bigr)\, =\, 
\lim_{n\to\infty} \frac{1}{c_{n\lambda}^*} P_{n\lambda}\bigl(\frac{z}{n}\bigr).\end{equation}
\end{remark}
 
For illustration, consider the rank-one case (type $A_1$) 
with $\f a =\b R$ and $R_+ =\{2\alpha\}, \,\alpha =1$. 
Fix $k=k_{2\alpha} \geq 0$. Then according to the example on p.89f of
\cite{O1},
\begin{align} F(\lambda,z) = \,&\,
   \phantom{}_2F_1 \bigl(a,b,c;\frac{1}{2}(1-\cosh\!z)\bigr), \notag\\
   G(\lambda,z) = \,&\,\phantom{}_2F_1
      \bigl(a,b,c;\frac{1}{2}(1-\cosh\!z)\bigr)\, + \,\frac{a}{2c}\sinh\!z\cdot 
      \phantom{}_2F_1 \bigl(a+1,b+1,c+1;\frac{1}{2}(1-\cosh\!z)\bigr),  \notag
\end{align}
with $\,a=\lambda+k, \,b=-\lambda+k,\, c = k+\frac{1}{2}.$ 
The weight lattice is $P=\b Z$, and the associated  Heckman-Opdam polynomials are given by
\begin{align} P_n(z) =\,&\, c_n^*F(n+k,z) =\, c_n^* \cdot Q_n^k(\cosh \!z), \,\, n = 0,1,\ldots \notag\\
E_n(z) =\,&\,c_n G(\widetilde n,z) = c_n\Bigl[Q_{|n|}^k(\cosh\!z) + \frac{\widetilde n +k}{2k+1} \cdot\sinh\!z\,Q_{|n|-1}^{k+1}(\cosh\!z)\Bigr],\,\, n\in \b Z \notag
\end{align}
with $\,\widetilde n = n+k$ for $n>0, \, \widetilde n = n-k$ for $n\leq 0$, and the renormalized Gegenbauer polynomials 
\[ Q_n^k(x) = \phantom{}_2F_1 \bigl(n+2k,-n,k+\frac{1}{2};\frac{1}{2}(1-x)\bigr).\] 
Relation \eqref{limit2} reduces to the classical limit 
\[ \lim_{n\to\infty} Q_n^k\bigl(\cos\frac{z}{n}\bigr) = 
j_{k-\frac{1}{2}}(z)
\quad(z\in \b C)\]
for the modified Bessel functions
\[ j_\alpha(z) = 2^\alpha\Gamma(\alpha+1)\cdot\frac{J_\alpha(z)}{z^\alpha}\,=\,
\Gamma(\alpha+1)\cdot\sum_{n=0}^\infty
\frac{(-1)^n(z/2)^{2n}}{n!\,\Gamma(n+\alpha+1)}\,,\]
see Theorem 4.11.6 of \cite{AAR}.
It is clear from the explicit  representation of the Gegenbauer polynomials in terms of Tchebycheff polynomials (Eq.~(6.4.11) in 
\cite{AAR}) that for $k\geq 0$, the expansion coefficients of $P_n$ with respect to the exponentials 
$z\mapsto e^{mz}, \, m\in \b Z$ are all non-negative. A closer inspection shows that the same holds for the non-symmetric $E_n$. 
This is in fact a special case
of a deep result for general Heckman-Opdam polynomials due to Sahi \cite{S}:
If the multiplicity function $k$ is non-negative, then it follows from Corollary 5.2. and Proposition 6.1. in \cite{S} that
the coefficients $a_{\lambda,\nu}$ of $E_\lambda$ in \eqref{coeff1} are all 
real and non-negative.
More precisely, if $\Pi_k := \b Z_+[k_\alpha]$ denotes the set 
of polynomials in the parameters $k_\alpha$ with non-negative integral 
coefficients, then
for suitable $d_\lambda\in \Pi_k$, all coefficients of 
$d_\lambda E_\lambda$ are contained in $\Pi_k$ as well. 
This positivity result is the key for our subsequent proof of Theorem \ref{postheorem}.

\section{New proof of Theorem \ref{postheorem}}

In contrast to our approach in \cite{R1}, we first derive
 a positive integral representation for the Dunkl kernel. As before, $R$
and $k\geq 0$ are fixed.

\begin{proposition}\label{kernelpos}
For each $x\in \f a$, there exists a unique probability measure $\mu_x\in M^1(\f a)$ such that
\begin{equation}\label{kernelrep}
 {\rm Exp}_W(x,z) = \int_{\f a} e^{\langle\xi,z\rangle} d\mu_x(\xi) 
\quad\forall z\in \f a_{\b C}.
\end{equation}
The support of $\mu_x$ is contained in $C(x)$. 
\end{proposition}

\begin{proof} It suffices to prove the existence of the representing measures as stated; their uniqueness  is immediate from the 
injectivity of the (usual) Fourier-Stieltjes transform on $M^1(\f a)$. 
Let $\lambda\in P$. Then by Sahi's positivity result mentioned above, 
\[ G(\widetilde\lambda,\,.\,) = \frac{1}{c_\lambda} E_\lambda\, =\, 
    \sum_{\nu\trianglelefteq \lambda} b_{\lambda,\nu}\, e^\nu  \]
with coefficients $b_{\lambda,\nu}$ satisfying
\[0\leq b_{\lambda,\nu} \leq 1, \quad\sum_{\nu\trianglelefteq \lambda}b_{\lambda,\nu} = 1.\]  
Now fix $\lambda\in P$ and $z\in \f a_{\b C}$. Then by Corollary \ref{limit},
\[ \text{Exp}_W(\lambda,z) = \, \lim_{n\to\infty} \frac{1}{c_{n\lambda}}E_{n\lambda}\bigl(\frac{z}{n}\bigr)\,=\, \lim_{n\to\infty} \sum_{\nu\trianglelefteq\, n\lambda} 
b_{n\lambda,\nu}\, e^{\langle\nu,\frac{z}{n}\rangle}.\]
Introducing the discrete probability measures
\[ \mu_\lambda^n := \sum_{\nu\trianglelefteq\, n\lambda}b_{n\lambda,\nu}\,
\delta_{\frac{\nu}{n}}\quad\in M^1(\f a)\] 
(where $\delta_x$ denotes the point measure in
  $x\in \f a$),
we may write the above relation in the form
\begin{equation}\label{int1}
  \text{Exp}_W(\lambda,z) = \, \lim_{n\to\infty}\int_{\f a}
  e^{\langle\xi,z\rangle} d\mu_\lambda^n(\xi).\end{equation}
The following Lemma shows that the support of $\mu_\lambda^n$ is contained 
in $C(\lambda).$

\begin{lemma}
Let $\lambda, \,\nu\in P$ with $\nu\trianglelefteq\, \lambda$. Then 
$\nu\in C(\lambda).$
\end{lemma}

\begin{proof}
Let $C:= \{x\in \f a: \langle\alpha,x\rangle\geq  0 \,\,\forall\, \alpha\in
R_+\}$
be the closed Weyl chamber associated with $R_+$ and
\[ C^* := \{y\in\f a: \langle y,x\rangle\geq 0 \quad\forall\, x\in C\}\]
its closed dual cone. Notice that $Q_+ \subset C^*$. Therefore
 $\nu\trianglelefteq\, \lambda$ implies that $\lambda_+  -\nu_+ \in C^*$.
We employ the following characterization of $C(x)$ for
$x\in C$ (Lemma IV.8.3 of \cite{He}):
\[ C(x) = \bigcup_{w\in W} w\bigl(C\cap(x-C^*)\bigr).\]
This shows that $\nu\in C(\lambda)$ if and only if
$\,\nu_+ \in \lambda_+ - C^*\,$, which yields the statement.
\end{proof}

We now continue with the proof of Proposition \ref{kernelpos}.
Fix $\lambda\in P$. By the preceding result, we may consider 
the  $\mu_\lambda^n$ as probability measures on the compact set $C(\lambda)$.
According to Prohorov's theorem (see e.g. \cite{Bi}), the set 
$\{\mu_\lambda^n, \, n\in \b Z_+\}$ is relatively compact. Passing 
 to a subsequence if necessary, we may therefore assume that 
there exists a measure $\mu_\lambda\in M^1(\f a)$ which is supported
in $C(\lambda)$ and such that $\mu_\lambda^n \to \mu_\lambda$ weakly as $n\to
 \infty$. Thus in view of \eqref{int1}, 
 \[\text{Exp}_W(\lambda,z) = \,\int_{\f a}
  e^{\langle\xi,z\rangle} d\mu_\lambda(\xi) \quad\forall 
  z\in \f a_{\b C}.\]
In order to extend this  representation  to arbitrary arguments $x\in \f a$ instead of $\lambda\in P$, observe first that 
for  $r\in \b Q$,
 \[\text{Exp}_W(r\lambda,z) =\text{Exp}_W(\lambda,rz) =\int_{\f a}
  e^{r\langle\xi,z\rangle} d\mu_\lambda(\xi).\]
Defining $\mu_{r\lambda}\in M^1(\f a)$ as the image measure of $\mu_\lambda$ under the dilation $\xi\mapsto r\xi$ on $\f a$, we therefore obtain \eqref{kernelrep}
for all $\,x\in \b Q.P 
=\{r\lambda: r\in \b Q\,, \lambda\in P\}.$
The set $\b Q.P$ is obviously dense in $\f a$. For arbitrary
 $x\in \f a$, choose an approximating sequence 
$\{x_n, \, n\in \b Z_+\} \subset \b Q.P$  with 
$\lim_{n\to\infty} x_n = x$.  Using Prohorov's theorem once more we obtain,
after passing to a subsequence, that $\mu_{x_n}\to \mu_x$ weakly for some
$\mu_x\in M^1(\f a)$. The support of $\mu_x$ can be confined to an arbitrarily small neighbourhood of $C(x)$, and must therefore coincide with $C(x)$. 
We thus have
\[ \text{Exp}_W(x,z) = \lim_{n\to\infty} \text{Exp}_W(x_n,z) = 
\int_{\f a}  e^{\langle\xi,z\rangle} d\mu_x(\xi) \quad \forall
z\in \f a_{\b C}, \]
which finishes the proof of the proposition.
\end{proof}

 \begin{proof}[Proof of Theorem \ref{postheorem}] 
By Proposition \ref{kernelpos} 
and the definition of $V$,
\[ \sum_{n=0}^\infty \frac{1}{n!} V_x\bigl(\langle x,z\rangle^n\bigr) =
V_x\bigl(e^{\langle x,z\rangle}\bigr) = 
   \int_{\f a}  e^{\langle\xi,z\rangle} d\mu_x(\xi) = \sum_{n=0}^\infty 
\frac{1}{n!}\int_{\f a} \langle\xi,z\rangle^n d\mu_x(\xi) \quad (z\in \f a_{\b C});\]
here the subscript $x$ means that $V$ is taken with respect to $x$.
Comparison of the homogeneous parts in $z$ of degree $n$ yields
that
\[ V_x(\langle x,z\rangle^n) = \int_{\f a} \langle\xi,z\rangle^n d\mu_x(\xi) \quad \forall n\in \b Z_+\,. \]
As the $\b C$-span of  $\{ x\mapsto \langle x,z\rangle^n, \, 
 z\in \f a_{\b C}\}$ is $\mathcal P_n$, it follows by linearity that 
\[ Vp(x) = \int_{\f a} p(\xi) d\mu_x(\xi) \quad\forall p\in \mathcal P, \, x\in \f a.\]
Finally, as $\mathcal P$ is dense in each $(A_r, \|\,.\,\|_{A_r})$ and 
$\|\,.\,\|_{\infty,K_r}\,\leq \|\,.\,\|_{A_r}$, an easy approximation
argument implies that this integral representation remains valid 
for all $f\in A_r$ with $r\geq |x|.$ This finishes the proof.
\end{proof}

We conclude this note with a remark concerning positive product formulas.
It is conjectured that (again in case $k\geq 0$) the multivariable 
hypergeometric function $F$ has a 
positive product formula. More precisely, we conjecture that
for all $x,y\in \f a$ there exists a probability measure $\sigma_{x,y}\in M^1(\f a)$ whose support is contained in the ball $K_{|x| + |y|}(0)$ and which satisfies
\begin{equation}\label{conj}
 F(\lambda,x)F(\lambda,y) = \int_{\f a} F(\lambda, \xi) d\sigma_{x,y}(\xi) 
\quad\forall \lambda\in \f a_{\b C}.\end{equation}
In the rank one case, i.e. for Jacobi functions, this is well-known and 
goes back to \cite{FK}. 
Eq.\eqref{conj} would immediately imply a positive product formula for the generalized Bessel function $J_W$ (associated with the same multiplicity $k$). In fact,
suppose there exist measures $\sigma_{x,y}$ as conjectured above, and 
denote for $r>0$  the image measure of  $\sigma_{x,y}$ under
the dilation $\xi\mapsto r\xi$ on $\f a$ by $\sigma^r_{x,y}$.
Then by  relation \eqref{limit2}, 
\begin{equation}\label{rephalf}
J_W(\lambda,x) J_W(\lambda,y) = \lim_{n\to\infty} 
   F\bigl(n\lambda + \rho, \frac{x}{n}\bigr)
   F\bigl(n\lambda +\rho , \frac{y}{n}\bigr) \,=\, \lim_{n\to\infty} 
\int_{\f a} F\bigl(n\lambda +\rho, \frac{\xi}{n}\bigr)d\sigma_{\frac{x}{n},\frac{y}{n}}^n(\xi)
\end{equation}
for all $\lambda\in \f a_{\b C}.$
As $\,\text{supp}\,\sigma_{\frac{x}{n},\frac{y}{n}}^n(\xi)\,\subseteq K_{|x|+|y|}(0)$ for all $n\in \b N$, we may assume that there exists a probability measure 
$\tau_{x,y}\in M^1(\f a)$ with $\text{supp}\,\tau_{x,y}\subseteq K_{|x|+|y|}(0)$ 
such that $\,\sigma_{\frac{x}{n},\frac{y}{n}}^n \to \tau_{x,y}$ weakly as $n\to\infty$. 
As further $\,\lim_{n\to\infty} F\bigl(n\lambda + \rho, \frac{\xi}{n}\bigr) = J_W(\lambda,\xi)$  locally  uniformly with respect to  $\xi$, 
Eq. \eqref{rephalf} implies the product
formula 
\[J_W(\lambda,x) J_W(\lambda,y) = \int_{\f a} J_W(\lambda,\xi) d\tau_{x,y}(\xi)\quad\forall \,\lambda\in \f a_{\b C}.\]
The uniqueness of $\tau_{x,y}$ is immediate from the injectivity 
of the Dunkl transform on $M^1(\f a)$, c.f. Theorem 2.6. of \cite{RV}.

\bigskip\medskip
\textbf{Acknowledgements.} We thank the referee for some useful
suggestions which helped to improve the presentation.
 The first author was supported by the Netherlands Organisation for Scientific Research (NWO), project number B 61-544.

\medskip

\noindent
Margit R\"osler\\
 Korteweg-de Vries Institute, University of Amsterdam\\
 Plantage Muidergracht 24\\
 1018 TV Amsterdam, The Netherlands\\
e-mail: mroesler@science.uva.nl

\medskip

\noindent
Michael Voit\\
Fachbereich Mathematik\\
Universit\"at Dortmund\\
Vogelpothsweg 87\\
44221 Dortmund, Germany\\
e-mail: michael.voit@mathematik.uni-dortmund.de

\end{document}